\input amstex



\def\b1{\text{\bf 1}}
\def\B1{{\Bbb 1}}

\def\BN{{\Bbb N}}

\def\BR{{\Bbb R}}
\def\BZ{{\Bbb Z}}

\def\CP{{\Cal P}}


\def\lra{\longrightarrow}
\def\ra{\rightarrow}

\parskip=6pt

\documentstyle{amsppt}
\document
\NoBlackBoxes



\centerline{NOMBRES DE BERNOULLI ET UNE FORMULE DE RAMANUJAN}

\bigskip\bigskip

\centerline{Oleg Ogievetsky\footnote{Centre de Physique Th\'eorique, Luminy, 
13288 Marseille (Unit\'e Mixte de Recherche 6207 du CNRS et des Universit\'es Aix--Marseille I, 
Aix--Marseille II et du Sud Toulon -- Var~; laboratoire affili\'e \`a la FRUMAM, FR 2291)} et 
Vadim Schechtman\footnote{Institut de Math\'ematique de 
Toulouse, Universit\'e Paul Sabatier, 31062 Toulouse}}

\vskip 1 cm

\centerline{\bf R\'esum\'e} 

\bigskip\bigskip

Dans la premi\`ere partie de cet article, on \'etablit une liaison \'etroite 
entre la formule de Euler - Maclaurin et l'\'equation 
fonctionelle 
de Rota - Baxter~; ces deux choses \'etant plus ou moins \'equivalentes. 

Dans la deuxi\`eme partie, on pr\'esente une simple d\'emonstration 
d'une formule de Ramanujan sur la sommation de certaines s\'eries 
exponentielles. 
Ceci a fait l'objet d'un expos\'e \`a Toulouse, en mai 2007. 

\bigskip\bigskip 

\centerline{ Table des mati\`eres} 

\bigskip\bigskip

Premi\`ere Partie. 

\'Equation de Rota - Baxter et formule sommatoire d'Euler - Maclaurin

\bigskip

\S 1. D\'efinition de Jacob Bernoulli . . . 2 

\S 2. Une primitive et l'\'equation de Rota - Baxter homog\`ene . . . 4

\S 3. Une primitive discr\`ete et l'\'equation de Rota - Baxter 
non-homog\`ene . . . 5

\S 4. Formule sommatoire d'Euler - Maclaurin . . . 7

Bibliographie . . . 10

\bigskip

Deuxi\`eme Partie. Une formule de Ramanujan

\bigskip

\S 1. Fonction $\eta$ de Dedekind . . . 11

\S 2. Une formule de Schl\"omilch . . . 14

\S 3. D\'eveloppements eul\'eriennes de $\sin$ et de $\cot$ . . . 15 

\S 4. Une formule de Ramanujan . . . 19 

\S 5. Une int\'egrale de Legendre . . . 22 

Bibliographie . . . 27

\bigskip\bigskip

\newpage

\centerline{Premi\`ere Partie}

\bigskip\bigskip

\centerline{\'EQUATION DE ROTA - BAXTER}

\bigskip 

\centerline{ ET FORMULE SOMMATOIRE D'EULER - MACLAURIN}

\bigskip\bigskip

\vskip 2 cm

\centerline{\bf \S 1. D\'efinition de Jacob Bernoulli} 

\bigskip\bigskip

{\bf 1.1.} Les nombres qu'A. de Moivre, puis Euler, ont appel\'es 
{\it nombres 
de Bernoulli}, ont \'et\'es introduits par Jacob I Bernoulli (1655 - 1705), 
dans son livre {\it Ars Conjectandi} sur les probabilit\'es,   
cf.  [B], Pars secunda, Caput III, pp. 97 - 98. Ce livre a \'et\'e publi\'e 
en 1713, quand Euler avait 6 ans (Euler fut un \'el\`eve du fr\`ere de
Jacob, Johann, et un ami de ses deux fils, Nicolas et Daniel).  

Bernoulli commence par un 
calcul de polyn\^omes qu'il d\'esigne par $\int n^r$~; nous convenons 
de la notation 
$$
S_r(n) = 1^r + 2^r + \ldots + n^r
$$
La m\'ethode de calcul est bas\'ee sur le triangle de Pascal 
(qui \`a l'\'epoque a servi pour la d\'efinition des {\it numerorum
figuratorum}, alias coefficients binomiaux). Cette m\'ethode \'etait 
d\'ej\`a connue de Pierre de Fermat. 

Voici ce qu'\'ecrit Bernoulli~:

\bigskip

"... Atque si porr\`o ad altiores gradatim potestates pergere, levique negotio 
sequentem adornare laterculum licet~: 

\bigskip 

\centerline{\it Summae Potestatum} 

$\int n = \frac{1}{2} nn + \frac{1}{2} n$

$\int nn = \frac{1}{3} n^3 + \frac{1}{2} nn + \frac{1}{6} n$ 

$\int n^3 = \frac{1}{4} n^4 + \frac{1}{2} n^3 + \frac{1}{4} nn$ 

$\int n^4 = \frac{1}{5} n^5 + \frac{1}{2} n^4 + \frac{1}{3} n^3 
- \frac{1}{30} n$

$\int n^5 = \frac{1}{6} n^6 + \frac{1}{2} n^5 + \frac{5}{12} n^4 
- \frac{1}{12} nn$ 

$\int n^6 = \frac{1}{7} n^7 + \frac{1}{2} n^6 + \frac{1}{2} n^5 
- \frac{1}{6} n^3 + \frac{1}{42} n$

$\int n^7 = \frac{1}{8} n^8 + \frac{1}{2} n^7 + \frac{7}{12} n^6 
- \frac{7}{24} n^4 + \frac{1}{12} nn$

$\int n^8 = \frac{1}{9} n^9 + \frac{1}{2} n^8 + \frac{2}{3} n^7 
- \frac{7}{15} n^5 + \frac{2}{9} n^3 - \frac{1}{30} n$

$\int n^9 = \frac{1}{10} n^{10} + \frac{1}{2} n^9 + \frac{3}{4} n^8 
- \frac{7}{10} n^6 + \frac{1}{2} n^4 - \frac{1}{12} nn$

$\int n^{10} = \frac{1}{11} n^{11} + \frac{1}{2} n^{10} + \frac{5}{6} n^9 
- 1\ n^7 + 1\ n^5 - \frac{1}{2} n^3 + \frac{5}{66} n$

Quin im\`o qui legem progressionis inibi attentuis ensperexit, eundem etiam 
continuare poterit absque his ratiociniorum ambabimus~: Sumt\^a enim $c$ 
pro potestatis cujuslibet exponente, fit summa omnium $n^c$ seu 
$$
\int n^c = \frac{1}{c+1} n^{c+1} + \frac{1}{2} n^{c} + 
\frac{c}{2}A n^{c-1} + \frac{c\cdot c-1\cdot c-2}{2\cdot 3\cdot 4}B n^{c-3}   
$$
$$
+ \frac{c\cdot c-1\cdot c-2\cdot c-3\cdot c-4}{2\cdot 3\cdot 4\cdot 5\cdot 6}
C n^{c-5}  
$$
$$
+ \frac{c\cdot c-1\cdot c-2\cdot c-3\cdot c-4\cdot c-5\cdot c-6}
{2\cdot 3\cdot 4\cdot 5\cdot 6\cdot 7\cdot 8}D n^{c-7} \ldots \&
\ \text{ita deinceps},
$$
exponentem potestatis ipsius $n$ continu\'e minuendo binario, quosque 
perveniatur ad $n$ vel $nn$. Literae capitales $A, B, C, D\ \&$ c. ordine
denotant co\"efficientes ultimorum terminorum pro $\int nn, \int n^4,   
\int n^6, \int n^8,\ \&$ c. nempe 
$$
A = \frac{1}{6},\ B = - \frac{1}{30},\ C = \frac{1}{42},\ D = - \frac{1}{30}. 
$$
Sunt autem hi coefficientes ita comparati, ut singuli cum caeteris sui 
ordinis co\"effi\-cien\-ti\-bus complere debeant unitatem~; sic $D$ valere 
diximus $-1/30$, 
$$
\text{quia\ } \frac{1}{9} + \frac{1}{2}  + \frac{2}{3}  
- \frac{7}{15}  + \frac{2}{9} (+D) - \frac{1}{30} = 1.
$$
Huius laterculi beneficio intra semi-quadrantem horae reperi, qu\`od potestates 
decime sive quadrato-sursolida mille primorum numerorum ab unitate in summam 
collecta efficiunt 
$$
91\ 409\ 924\ 241\ 424\ 243\ 424\ 241\ 924\ 242\ 500."
$$

\bigskip 

{\bf 1.2.} D\'efinissons les nombres $b_n$ par la s\'erie g\'en\'eratrice 
$$
\frac{S}{1 - e^{-S}} = \sum_{n=0}^\infty\ 
\frac{b_n}{n!}S^n = 1 + \frac{S}{2} + \sum_{p=1}^\infty\ 
\frac{b_{2p}}{(2p)!}S^{2p}\ .
$$
On remarque que 
$$
\frac{S}{e^S - 1} = 1 - \frac{S}{2} + \sum_{p=1}^\infty\ 
\frac{b_{2p}}{(2p)!}S^{2p}\ .
$$
Voici les premiers valeurs~: 
$$
b_2 = \frac{1}{6},\  b_4 = - \frac{1}{30},\  
b_6 = \frac{1}{42},\  b_8 = - \frac{1}{30},\ 
$$
$$
b_{10} = \frac{5}{66},\  b_{12} = - \frac{691}{2730},\  
b_{14} = \frac{7}{6}\ .
$$

{\bf 1.3.} On verra plus bas que $S_r(n)$ est la valeur en $x=n$ du polyn\^ome 
$$
S_r(x) = \frac{x^{r+1}}{r+1} + \frac{x^r}{2} + \sum_{1\leq p < (r+1)/2}\ 
\binom{r}{2p-1}\frac{b_{2p}}{2p}x^{r-2p+1}\ .
$$
Autrement dit, 
$$
S'_r(x) = B_r(x) = \sum_{p=0}^r\ \binom{r}{p}b_px^{r-p},\ \ \ S_r(0) = 0\ .
$$

\bigskip\bigskip

\centerline{\bf \S 2. Une primitive et l'\'equation de Rota-Baxter homog\`ene} 

\bigskip\bigskip 

{\bf 2.1.} Soit $A$ une alg\`ebre de fonctions $f(x)$ "raisonnables", 
par exemple l'alg\`ebre des polyn\^omes $\BR[x]$ ou l'alg\`ebre des fonctions 
d\'erivables. On d\'esigne par $D$ l'op\'e\-ra\-teur de d\'erivation sur $A$, 
et par $I$ l'op\'erateur 
$$
I(f)(x) = \int_0^x f(t)dt \ .
$$
Il est clair que $DI =\ $id. En revanche, 
$$
ID(f)(x) = f(x) - f(0)\ .
$$

{\bf 2.2.} {\it Lemme.} L'op\'erateur $I$ satisfait l'\'equation 
$$
I(f)I(g) = I(I(f)g + fI(g))\ .
\eqno{(RBH)}
$$

{\it Premi\`ere preuve.} Les d\'eriv\'ees des deux c\^ot\'es co\^\i ncident 
puisque 
$DI =\,$ id. De plus, les valeurs des deux c\^ot\'es en $0$ sont $0$, 
d'o\`u l'assertion. 

{\it Seconde preuve.} Consid\'erez l'int\'egrale de la fonction 
de deux variables $f(t)g(u)$ sur le carr\'e $[0,x]^2$~; puis coupez ce carr\'e 
en deux triangles.  

L'\'equation (RBH) sera appel\'ee {\it \'equation de Rota-Baxter homog\`ene}, cf. [Ro].

\bigskip\bigskip

\newpage

\centerline{\bf \S 3. Une primitive discr\`ete et l'\'equation de Rota - Baxter 
non-homog\`ene} 

\bigskip\bigskip 

{\bf 3.1.} \'Etant donn\'ee une fonction $f: \BN \lra \BR$ d\'efinissons 
sa "primitive discr\`ete" $B(f):\ \BN \lra \BR$ par 
$$
B(f)(n) = \sum_{i=1}^n\ f(i) \ .
$$

{\bf 3.2.} {\it Lemme.} L'op\'erateur $B$ satisfait l'\'equation 
$$
B(fg) + B(f)B(g) = B(B(f)g + fB(g))\ .
\eqno{(RB)} 
$$

En effet, la valeur du premier membre de (RB) en $n$ est 
$$
\sum_{i=1}^n f(i)g(i) + \sum_{i=1}^n f(i)\cdot \sum_{j=1}^n g(j)\ .
$$
Donc c'est un carr\'e $n\times n$, avec la diagonale doubl\'ee. 
D'autre part, 
le deuxi\`eme membre est 
$$
\sum_{i=1}^n\biggl\{\sum_{j=1}^if(j)g(i) + \sum_{j=1}^if(i)g(j)\biggr\} \ .
$$
Il est ais\'e de voir que les deux expressions sont \'egales, en interpretant 
chaque terme $\sum_{j=1}^if(j)g(i) + \sum_{j=1}^if(i)g(j)$ comme un chemin 
(de forme $\Gamma$) dans le m\^eme carr\'e.  

{\bf 3.3.} Il est clair qu'on peut consid\'erer l'anneau de polyn\^omes 
$\BR[x]$ comme un sous-anneau de l'anneau $\BN^\BR$ des applications 
$f: \BN \lra \BR$. 

{\it Lemme} (Bernoulli). $B(\BR[x]) \subset \BR[x]$. 

En effet, on peut calculer les polyn\^omes $S_r(x) := B(x^r)$ par 
r\'ecurrence, en utilisant (RB)~: 

On a $B(1)(n) = n$, donc $B(1) = x$. Ensuite, 
$$
B(1\cdot 1) + B(1)B(1) = B(B(1)1 + 1B(1)), 
$$
i.e. 
$$
x + x^2 = 2B(x), 
$$
d'o\`u $B(x) = (x^2 + x)/2$. 

Ensuite, 
$$
B(1\cdot x) + B(1)B(x) = B(B(1)x + 1B(x)), 
$$
i.e. 
$$
\frac{x^2 + x}{2} + \frac{x^3 + x^2}{2} = B(x^2 + \frac{x^2 + x}{2}) = 
\frac{1}{2}B(3x^2 + x) = \frac{3}{2}B(x^2) + \frac{x^2 + x}{4}, 
$$
d'o\`u
$$
B(x^2) = \frac{x^3}{3} + \frac{x^2}{2} + \frac{x}{6}, 
$$
cf. 1.1. Ainsi, les $B(x^i)$ pour $i\leq r$ \'etant connus, on obtient 
$B(x^{r+1})$ en appliquant (RB) avec $f = 1, g = x^r$, ce qui prouve 
le lemme. 

\bigskip

{\bf 3.4.} Soit $A$ une alg\`ebre associative munie d'un op\'erateur lin\'eaire
$B: A \lra A$ v\'erifiant 
$$
\mu B(fg) + B(f)B(g) = B(B(f)g + fB(g))\ 
\eqno{(RB)}
$$
($\mu$ est un nombre).
On introduit sur $B$ une autre multiplication 
$$
f * g = B(f)g + fB(g) -\mu fg\ .
$$
En l'utilisant, on peut r\'eecrire (RB) sous la forme \'equivalente 
$$
B(f*g)  = B(f)B(g)\ .
\eqno{(RB)'}
$$
Autrement dit, $B$ est un op\'erateur entrela\c{c}ant deux multiplications.  

{\bf 3.5.} {\it Lemme.} La multiplication $*$ est associative. 

Se v\'erifie aisement \`a l'aide de $(RB)'$. 

{\bf 3.6.} {\it Lemme.} L'op\'erateur $B$ satisfiait \`a $(RB)$ pour $*$~: 
$$
\mu B(f*g) + B(f)*B(g) = B(B(f)*g + f*B(g))\ .
$$

En effet, 
$$
B(f*g) + B(f)*B(g) = B^2(f)B(g) + B(f)B^2(g)\ .
$$
D'un autre c\^ot\'e, d'apr\`es $(RB)'$, on a  
$B(B(f)*g)) = B^2(f)B(g)$ et $B(f*B(g)) = B(f)B^2(g)$, d'o\`u l'assertion.  

\bigskip\bigskip

\newpage

\centerline{\bf \S 4. Formule sommatoire d'Euler - Maclaurin} 

\bigskip\bigskip 

{\bf 4.1.} Soit $A$ une alg\`ebre commutative munie d'une d\'erivation $D$ 
et d'un op\'erateur $I:\ A \lra A$ ("une primitive") tel que $DI =\ $id$_A$ et 
satisfaisant 
$$
I(f)I(g) = I(I(f)g + fI(g)) \ .
\eqno{(RBH)} 
$$
On veut construire un op\'erateur $B:\ A\lra A$ ("une primitive discr\`ete") 
de la forme 
$$
B = I(1 + \beta_1 D + \beta_2 D^2 + \ldots),\ \beta_i \in \BR\ ,
\eqno{(4.1.1)}
$$
qui satisfait 
$$
B(fg) + B(f)B(g) = B(B(f)g + fB(g))\ .
\eqno{(RB)} 
$$
Essayons de trouver un par un les coefficients inconnus $\beta_i$, 
de l'\'equation (RB). 

{\bf 4.2.} On a \`a gauche~: 
$$
B(fg) = I(fg) + \beta_1 ID(fg) + \beta_2 ID^2(fg) + \ldots 
$$
$$
= I(fg) + \beta_1 I(Df\cdot g + f\cdot Dg) + 
\beta_2 I(D^2f\cdot g + 2Df\cdot Dg + f\cdot D^2g) + \ldots
$$
Ensuite,
$$
B(f)B(g) = (I(f) + \beta_1ID(f) + \beta_2ID^2(f) + \ldots)\cdot 
(I(g) + \beta_1ID(g) + \beta_2ID^2(g) + \ldots) 
$$
$$
= I(f)I(g) + \beta_1(ID(f)I(g) + I(f)ID(g)) 
$$
$$ 
+ \beta_2(ID^2(f)I(g) + I(f)ID^2(g)) + 
\beta_1^2 ID(f)ID(g) + \ldots 
$$
$$
= I\{I(f)g + fI(g)\} 
$$    
$$
 + \beta_1 I\bigl\{ID(f)g + Df I(g) + I(f)Dg + fID(g)\bigr\} 
$$
$$
+ \beta_2 I\bigl\{ID^2(f)g + D^2(f)I(g) + I(f)D^2g + fID^2(g)\} 
$$
$$
+ \beta_1^2 I\bigl\{ID(f)D(g) + D(f)ID(g)\} + \ldots 
$$

{\bf 4.3.} \`A droite~: 
$$
B(B(f)g + fB(g)) 
$$
$$
= I \biggl\{I(f)g + fI(g) + \beta_1(ID(f)g + f\cdot ID(g)) + 
\beta_2(ID^2(f)g + f\cdot ID^2(g)) + \ldots\biggr\} 
$$
$$
+ \beta_1 ID\biggl\{
I(f)g + fI(g) + \beta_1(ID(f)g + f\cdot ID(g)) + \ldots\biggr\} 
$$
$$
+ \beta_2 ID^2\biggl\{I(f)g + fI(g) + \ldots\biggr\} + \ldots 
$$
$$
= I \biggl\{I(f)g + fI(g) + \beta_1(ID(f)g + f\cdot ID(g)) + 
\beta_2(ID^2(f)g + f\cdot ID^2(g)) + \ldots\biggr\} 
$$
$$
+ \beta_1 I \biggl\{fg + I(f)D(g) + D(f)I(g) + fg\biggr\} 
$$
$$
+ \beta_1^2I\biggl\{D(f)g + ID(f)D(g) + D(f)ID(g) + fD(g)\biggr\} 
$$
$$
+ \beta_2 I\biggl\{D(f)g + 2f D(g) + I(f)D^2(g) + 
D^2(f)I(g) + 2D(f)g + fD(g)\biggr\} + \ldots
$$
(on a gard\'e les termes d'ordre $\leq 2$ en $D$). 

{\bf 4.4.} En r\'ealisant l'\'equation, les seuls termes qui surviennent sont~: 
$I(fg)$, ce qui donne 
$$
1 = 2\beta_1,
\eqno{(4.4.1)}
$$ i.e. $\beta_1 = 1/2 = - b_1$~; puis,  

$I(D(f)g + fD(g))$, qui donne 
$$
\beta_1 = 3\beta_2 + \beta_1^2,
\eqno{(4.4.2)}
$$ 
d'o\`u 
$$
\beta_2 = \frac{1}{12} = \frac{1}{2\cdot 6}\ .
$$ 
De la m\^eme mani\`ere, les calculs \`a l'ordre $3$ fournissent la valeur 
$\beta_3 = 0$. Ils donnent \`a l'ordre $4$, en regardant les termes 
$I(D^3(f)g + fD^3(g))$ (et aussi $I(2D^2(f)D(g) + 2D(f)D^2(g))$), la r\'elation 
$$
0 = 5\beta_4 + \beta_2^2,
\eqno{(4.4.3)}
$$
d'o\`u
$$
\beta_4 = - \frac{\beta_2^2}{5} = - \frac{1}{5\cdot 144} = 
- \frac{1}{4!\cdot 30} \ .
$$
Ceci nous am\`ene \`a l'id\'ee que 
$$
\beta_n = \frac{b_n}{n!}, \ n \geq 2.  
$$
Autrement dit, on attend que la s\'erie g\'en\'eratrice des nombres 
$\beta_n$ soit 
$$
\frac{S}{1 - e^{-S}} = \sum_{n=0}^\infty\ \beta_nS^n\ .
$$

\bigskip

{\bf 4.5.} \'Ecrivons la s\'erie (4.1.1) comme 
$$
B = I\mu, 
$$
o\`u $\mu = \mu(D) \in \BR[[D]]$, $\mu(0) = 1$. L'\'equation (RB)  
s'\'ecrit alors~: 
$$
I\mu(f)I\mu(g) + I\mu(fg) = I\mu(I\mu(f)\cdot g + f\cdot I\mu(g))\ .
$$
On pose~: $\mu(f) = a,\ \mu(g) = b, \nu = \mu^{-1}$, donc $f = \nu(a), 
g = \nu(b)$. On obtient~: 
$$
Ia \cdot Ib + I\mu(\nu(a)\nu(b)) = I\mu(Ia\cdot \nu(b) + \nu(a)Ib)\ .
$$
On multiplie les deux c\^ot\'es par $\nu D$ (en prenant en compte le fait que $DI =\ $id)~: 
$$
\nu(a Ib + Ia\cdot b) + \nu(a)\nu(b) = Ia\cdot \nu(b) + \nu(a)Ib\ .
$$ 
Ensuite, on pose $Ia = \alpha,\ Ib = \beta$, donc $a = D\alpha,\ b = D\beta$~: 
$$
\nu D(\alpha\beta) + 
\nu D(\alpha)\nu D(\beta) = \alpha \nu D(\beta) + \nu D(\alpha) \beta\ .
$$
Il s'en suit que si l'on pose $\kappa = 1 - \nu D$, alors 
$$
\kappa(\alpha\beta) = \kappa(\alpha)\kappa(\beta)\ .
\eqno{(4.5.1)}
$$
Il est clair que $\kappa(1) = 1$. 

{\bf 4.6.} Maintenant supposons que notre alg\`ebre $A$ est l'anneau 
de polyn\^omes $A = \BR[x]$. L'identit\'e (4.5.1) implique que 
$\kappa$ est un automorphisme de $A$, donc il est de la forme 
$\kappa(x) = Ex + F$. 

De plus, on impose une condition de normalisation 
$$
B(1) = x,
\eqno{(4.6.1)} 
$$
d'o\`u $I\mu(1) = x$, donc $\mu(1) = 1$, donc $\nu(1) = 1$, d'o\`u 
$$
\kappa(x) = (1 - \nu D)(x) = x - 1\ .
$$
Il vient que $\kappa = e^{-D}$, d'o\`u $\nu D = 1 - e^{-D}$, 
$\nu = (1 - e^{-D})D^{-1}$, donc 
$$
\mu(D) = \frac{D}{1 - e^{-D}}, 
$$
comme attendue. Donc 
$$
B = \frac{ID}{1 - e^{-D}}\ .
\eqno{(4.6.2)}
$$

{\bf 4.7.} Explicitement, on a~: 
$$
\frac{ID}{1 - e^{-D}} = I + \frac{ID}{2} + \sum_{p=1}^\infty\ 
\frac{b_{2p}}{(2p)!} ID^{2p}\ .
$$
Maintenant, pour $f(x) \in \BR[x]$ on a 
$$
ID(f)(x) = f(x) - f(0), 
$$
donc 
$$
\frac{ID}{1 - e^{-D}}(f)(x) = \int_0^n f(t)dt + 
\frac{1}{2}(f(x) - f(0)) 
+ \sum_{p=1}^\infty\ \frac{b_{2p}}{(2p)!}(f^{(2p-1)}(x) - f^{(2p-1)}(0)) \ .
$$
On a montr\'e que cet op\'erateur v\'erifie (RB) et (4.6.1), or, il n'existe 
qu'un seul op\'erateur de la sorte (l'unicit\'e se voit tout de suite par r\'ecurrence),
celui qui \`a $f\in \BR[x]$ fait correspondre $B(f)\in \BR[x]$ tel que 
$B(f)(n) = \sum_{i=1}^n f(i)$. 

Il vient le 

{\bf 4.8.} {\it Th\'eor\`eme} (Euler - Maclaurin). Pour tout $f\in \BR[x]$ et 
$n\in \BN$ on a 
$$
f(1) + f(2) + \ldots + f(n) 
$$
$$
= \int_0^n f(t)dt + 
\frac{1}{2}(f(n) - f(0)) 
+ \sum_{p=1}^\infty\ \frac{b_{2p}}{(2p)!}(f^{(2p-1)}(n) - f^{(2p-1)}(0))\ . 
$$

En l'appliquant \`a $f(x) = x^r$, on obtient 
$$
1^r + 2^r + \ldots + n^r = S_r(n), 
$$
o\`u le polyn\^ome $S_r(x)$ est d\'efini par
$$
S_r(x) = \frac{x^{r+1}}{r+1} + \frac{x^r}{2} + \sum_{1\leq p < (r+1)/2}\ 
\binom{r}{2p-1}\frac{b_{2p}}{2p}x^{r-2p+1}\ .
$$ 

{\bf 4.9.} Appliquons (RB) \`a $f = 1,\ g = x^r$~: 
$$
B(x^r) + xB(x^r) = B(x^{r+1} + B(x^r)), 
$$
i.e.
$$
B(x^{r+1}) = (x+1)B(x^r) - B(B(x^r))\ .
$$
Cette identit\'e est \'equivalente \`a une identit\'e classique pour 
les nombres de Bernoulli~: 
$$
(2n+1)b_{2n} = - \sum_{p=1}^{n-1}\ \binom{2n}{2p}b_{2p}b_{2n-2p},
$$ 
cf. [Bo], Ch. VI, \S 2, Exercice 2)~; [R] (c'est le premier article publi\'e de 
Ramanujan).    

\bigskip\bigskip

\bigskip\bigskip

\centerline{\bf Bibliographie}

\bigskip\bigskip

[B] Jacob Bernoulli, Ars conjectandi. B\^ale, 1713. (Dans~: Die Werke von Jacob Bernoulli,
Band 3, Birkh\"auser Verlag Basel, 1975.) 

[Bo] N.Bourbaki, Fonctions d'une variable r\'eelle. Hermann, Paris, 1961. 

[R] S.Ramanujan, Some properties of Bernoulli's numbers, {\it J. Ind. Math. 
Soc.}, III, 1911, 219 - 234. (Dans~: Collected papers, AMS Chelsea, 2000, pp. 1 - 14.) 
 
[Ro] G.-C. Rota, Baxter operators, an introduction. Dans~: Gian-Carlo Rota on combinatorics, 
Contemp. Mathematicians, Birkh\"auser, 1995, pp. 504-512.

\bigskip

\newpage

\centerline{Deuxi\`eme Partie} 

\bigskip\bigskip

\centerline{UNE FORMULE DE RAMANUJAN}


\vskip 3 cm

\centerline{\bf \S 1. Fonction $\eta$ de Dedekind} 

\bigskip\bigskip 

{\bf 1.1.} Il semble que Riemann ait lu assez attentivement les {\it Fundamenta} 
de Jacobi.   
Dans les papiers de Riemann, on a trouv\'e un {\it Additamentum ad} 
\S$^{\text{um}}$ 40 de {\it "Fundamenta"}, [R]. Richard Dedekind 
a \'ecrit 
un commentaire sur ces fragments, [D], o\`u il introduit 
la fonction $\eta(\tau)$~:   
$$
\eta(\tau) = q^{1/24}\prod_{n=1}^\infty\ (1 - q^n),\ 
q = e^{2\pi i\tau},
\eqno{(1.1.1)} 
$$
o\`u $|q| < 1$, i.e. $\Im \tau > 0$, 
et \'etudie sa loi de transformation
par rapport aux transformations 
de Moebius $\tau \mapsto (a\tau + b)/(c\tau + d)$. Le th\'eor\`eme suivant 
en est un cas particulier.   
 
{\bf 1.2.} {\it Th\'eor\`eme.} La fonction $\eta(\tau)$ satisfait 
l'\'equation 
$$
\eta(-1/\tau) = \sqrt{\tau/i}\eta(\tau)\ .
\eqno{(1.2.1)} 
$$

{\it D\'emonstration}, d'apr\`es Carl Ludwig Siegel, [S]. En prenant 
le logarithme naturel, 
$$
\frac{\pi i\tau}{12} - \log\eta(\tau) = - \sum_{n=1}^\infty\ 
\log(1 - q^n) = \sum_{n,m = 1}^\infty\ \frac{q^{nm}}{m} = 
\sum_{m=1}^\infty\ \frac{1}{m(q^{-m} - 1)}\ .
$$
Prenons le logarithme de (1.2.1)~: 
$$
\log\eta(-1/\tau) = \frac{1}{2}\log(\tau/i)  + \log\eta(\tau), 
$$
ou
$$
\frac{\pi i\tau}{12} - \log\eta(\tau) = - \frac{\pi i \tau^{-1}}{12} 
- \log\eta(-1/\tau) + \frac{1}{2}\log(\tau/i) + 
\frac{\pi i(\tau + \tau^{-1})}{12}\ .
$$
Donc (1.2.1) est \'equivalente \`a~: 
$$
\frac{1}{2}\log(\tau/i) + 
\frac{\pi i(\tau + \tau^{-1})}{12} = 
\sum_{m=1}^\infty\ \frac{1}{m}\biggl(\frac{1}{e^{-2\pi im\tau} - 1} 
- \frac{1}{e^{2\pi im/\tau} - 1}\biggr)\ .
\eqno{(1.2.2)}
$$

{\bf 1.3.} {\it Une fonction int\'eressante~:} $\cot z$. On pose $y = e^{iz}$. 
Alors~: 
$$
\cot z = \frac{\cos z}{\sin z} = 
- \frac{(y^{-1} + y)/2}{(y^{-1} - y)/2i} = 
- i\cdot \frac{y^{-1} + y}{y^{-1} - y} = 
i\cdot \frac{y + y^{-1}}{y - y^{-1}} 
$$  
$$
= - i\cdot\biggl(1 + \frac{2}{y^{-2} - 1}\biggr) = 
i\cdot\biggl(1 + \frac{2}{y^{2} - 1}\biggr)\ .
\eqno{(1.3.1)}
$$
Donc $\lim_{y\ra 0}\cot z = - i$ et $\lim_{y\ra \infty}\cot z = i$. 
De l\`a~: 
$$
\lim_{n\ra\infty}\cot((n+1/2)z) = - i\ \text{si}\ \Im z > 0
\eqno{(1.3.2a)} 
$$
et
$$
\lim_{n\ra\infty}\cot((n+1/2)z) = i\ \text{si}\ \Im z < 0\ .
\eqno{(1.3.2b)} 
$$

{\bf 1.4.} On pose $f(z) = \cot z \cot z/\tau$ et on consid\`ere la 
fonction $g_n(z) = z^{-1}f(\nu z)$ o\`u $\nu = (n + 1/2)\pi$, 
$n = 0, 1, \ldots $ Soit $C$ le contour du parallelogramme de sommets 
$1, \tau, - 1, - \tau$. 

Quels sont les p\^oles de $g_n(z)$? On a~: 
$$
g_n(z) = \frac{\cos \nu z}{z\sin \nu z}\cdot 
\frac{\cos \nu z/\tau}{\sin \nu z/\tau}\ .
$$
Donc on a~: 

(a) des p\^oles simples en $z = \pm \pi m/\nu$, $m = 1, 2, \ldots$, 
avec les r\'esidus
$$
\text{res}_{z = \pm\pi m/\nu}\ g_n(z) = \frac{\cot(\pi m/\tau)}{\pi m}\ ; 
$$

(b) des p\^oles simples en $z = \pm \pi m\tau/\nu$, $m = 1, 2, \ldots$, 
avec les r\'esidus
$$
\text{res}_{z = \pm\pi m\tau/\nu}\ g_n(z) = \frac{\cot(\pi m\tau)}{\pi m}\ . 
$$

(c) Enfin, en $z=0$ on a~: 
$$
g_n(z) = \frac{1}{z}\cdot \frac{1}{\nu z}\cdot \frac{\tau}{\nu z}\cdot 
\frac{1 - \nu^2z^2/2 + \ldots}{1 - \nu^2z^2/6 + \ldots}\cdot 
\frac{1 - \nu^2z^2/2\tau^2 + \ldots}{1 - \nu^2z^2/6\tau^2 + \ldots} 
$$ 
$$
= \frac{\tau}{\nu^2z^3}\cdot\biggl(1 - \frac{\nu^2z^2}{3} + \ldots\biggr)\cdot 
\biggl(1 - \frac{\nu^2z^2}{3\tau^2} + \ldots\biggr) 
$$
$$
= \frac{\tau}{\nu^2z^3}\cdot\biggl(1 - \frac{\nu^2z^2}{3}\cdot 
(1 + \tau^{-2}) + \ldots \biggr), 
$$
d'o\`u
$$
\text{res}_{z = 0}\ g_n(z) = - \frac{\tau + \tau^{-1}}{3}\ .
$$
Par la formule des r\'esidus de Cauchy, 
$$
\frac{1}{2\pi i}\int_C\ f(\nu z)\frac{dz}{z} = - \frac{\tau + \tau^{-1}}{3} 
+ \frac{2}{\pi}\sum_{m=1}^n\ \frac{1}{m}(\cot\pi m\tau + \cot\pi m/\tau)\ .
$$
On remarque que 
$$
\cot\pi m\tau + \cot\pi m/\tau = - 2i\biggl(\frac{1}{e^{-2\pi im\tau} - 1}
- \frac{1}{e^{2\pi im/\tau} - 1}\biggr), 
$$
cf. (1.3.1), d'o\`u
$$
\int_C\ f(\nu z)\frac{dz}{z} = - \frac{2\pi i(\tau + \tau^{-1})}{3} 
+ 8\sum_{m=1}^n\ \frac{1}{m}\biggl(\frac{1}{e^{-2\pi im\tau} - 1}
- \frac{1}{e^{2\pi im/\tau} - 1}\biggr)\ .
\eqno{(1.4.1)}
$$

{\bf 1.5.} Maintenant faisons tendre $n$ \`a l'infini dans (1.4.1). 
Soit $\ell_1 = \{\Im z = 0\}$ et $\ell_2$ la droite qui passe par
$0$ et $\tau$. D'apr\`es (1.3.2a,b), 

$\lim_{n\ra\infty}\cot \nu z = - i$ si $z$ est au-dessus de $\ell_1$~;  
$\lim_{n\ra\infty}\cot \nu z =  i$ si $z$ est au-dessous de $\ell_1$ et

$\lim_{n\ra\infty}\cot \nu z/\tau = i$ si $z$ est \`a droite de $\ell_2$~;  
$\lim_{n\ra\infty}\cot \nu z/\tau = - i$ si $z$ est \`a gauche de $\ell_2$. 

Il s'en suit que sur le c\^ot\'e $(1,\tau)$ de $C$ (sans les sommets) 
la valeur limite $\lim_{n\ra\infty} \cot\nu z\cot\nu z/\tau = -i\cdot i = 1$. 

De m\^eme, sur les c\^ot\'es $(\tau,-1)$, $(-1,-\tau)$ et $(-\tau,1)$ 
les valeurs limites sont $-1, 1, -1$. 

De l\`a, 
$$
\lim_{n\ra\infty}\int_C\ f(\nu z)\frac{dz}{z} = 
\biggl(\int_1^\tau - \int_\tau^{-1} + \int_{-1}^{-\tau} - \int_{-\tau}^{1}
\biggr)\frac{dz}{z} 
$$
$$
= \log\tau  - \pi + \log\tau 
+ \log(-\tau) - \pi - 2\pi + \log(-\tau) = 4\log\tau - 2\pi 
= 4\log(\tau/i)\ .
\eqno{(1.5.1)}
$$
Donc en passant \`a la limite $n\ra\infty$ dans (1.4.1), on obtient~: 
$$
4\log(\tau/i) + \frac{2\pi i(\tau + \tau^{-1})}{3} = 
8\sum_{m=1}^n\ \frac{1}{m}\biggl(\frac{1}{e^{-2\pi im\tau} - 1}
- \frac{1}{e^{2\pi im/\tau} - 1}\biggr)\ .
$$
En divisant par $8$, on obtient la formule cherch\'ee (1.2.2), 
QED.

\bigskip\bigskip

\newpage

\centerline{\bf \S 2. Une formule de Schl\"omilch} 

\bigskip\bigskip

{\bf 2.1.} {\it Th\'eor\`eme}, [Sch], [Ram]. 
$$
\sum_{n=1}^\infty\ \frac{n}{e^{2\pi n} - 1} = \frac{1}{24} - \frac{1}{8\pi}\ .
\eqno{(2.1.1)} 
$$

{\bf 2.2.} {\it D\'emonstration} de Srinivasa Ramanujan, [Ram], (18), p. 32. 
On prend $\tau = ia$ dans (1.2.1), o\`u $a$ est un nombre r\'eel, $a > 0$~: 
$$
e^{-\pi/12a}\prod_{n=1}^\infty\ (1 - e^{-2\pi n/a}) = 
\sqrt a\cdot e^{-\pi a/12}\prod_{n=1}^\infty\ (1 - e^{-2\pi na})\ .
$$
En prenant le logarithme, 
$$
- \frac{\pi}{12a} + \sum_{n=1}^\infty\ \log(1 - e^{-2\pi n/a}) 
= \frac{\log a}{2} - \frac{\pi a}{12} + 
\sum_{n=1}^\infty\ \log(1 - e^{-2\pi na})\ .
$$
En prenant la d\'eriv\'ee, 
$$
\frac{\pi}{12a^2} - \sum_{n=1}^\infty\ 
\frac{(2\pi n/a^2)\cdot e^{-2\pi n/a}}{1 - e^{-2\pi n/a}} = 
\frac{1}{2a} - \frac{\pi}{12} + \sum_{n=1}^\infty\ 
\frac{2\pi ne^{-2\pi n/a}}{1 - e^{-2\pi na}}, 
$$
ou bien
$$
\frac{\pi}{12}(a^{-2} + 1) - \frac{1}{2a} = 2\pi \sum_{n=1}^\infty\  
n\cdot\biggl(\frac{a^{-2}}{e^{2\pi n/a} - 1} + 
\frac{1}{e^{2\pi na} - 1}\biggr)\ .
\eqno{(2.2.1)} 
$$
Sous une forme plus sym\'etrique,
$$
 \frac{\pi(a^{-1} + a)}{12} - \frac{1}{2} = 2\pi \sum_{n=1}^\infty\  
 \biggl(\frac{n/a}{e^{2\pi n/a} - 1} + 
 \frac{na}{e^{2\pi na} - 1}\biggr)\ .
 \eqno{(2.2.2)} 
$$

En posant $a=1$, on arrive \`a (2.1.1).

\bigskip\bigskip 

\newpage

\centerline{\bf \S 3. D\'eveloppements eul\'eriens de $\sin$ et de 
$\cot$} 

\bigskip\bigskip

{\bf 3.1.} On suit Bourbaki, [B], Chapitre VI, \S 2. 

{\it Lemme.} On a pour $n\in \BZ$, $n>0$~:
$$
\sin nz = 2^{n-1}\prod_{k=0}^{n-1}\ \sin(z+k\pi/n)\ .
$$

En effet, 
$$
\sin nz = \frac{1}{2i}(e^{inz} - e^{-inz}) = \frac{e^{-inz}}{2i} 
(e^{2inz} - 1) 
$$
$$
= \frac{e^{-inz}}{2i}\prod_{p=0}^{n-1}(e^{2iz} - e^{-2\pi ip/n}) = 
\frac{1}{2i}\prod_{p=0}^{n-1}(e^{iz} - e^{-iz-2\pi ip/n}) 
$$
$$
= (2i)^{n-1}\prod_{p=0}^{n-1}e^{-\pi ip/n}
\prod_{p=0}^{n-1}\frac{e^{iz+\pi ip/n} - e^{-iz-\pi ip/n}}{2i}\ .
$$
Or,
$$
(2i)^{n-1}\prod_{p=0}^{n-1}e^{-\pi ip/n} = (2i)^{n-1} 
e^{-\pi i/n\cdot \sum_{p=0}^{n-1} p} = (2i)^{n-1}e^{-\pi i(n-1)/2} = 
2^{n-1},
$$ 
d'o\`u l'assertion. 

{\bf 3.2.} En divisant par $\sin z$ et en faisant tendre $z$ vers $0$, on
obtient 
$$
\prod_{p=0}^{n-1}\ \sin(p\pi/n) = n 2^{1-n}\ .
$$

{\bf 3.3.} Soit $n = 2m + 1$ impair. On a~: 
$\sin(n(z + \pi/2)) = \sin(nz + \pi/2 + m\pi) = (-1)^m\cos nz$, d'o\`u, 
en rempla\c{c}ant $z$ par $z + \pi/2$ dans 3.1, 
$$
\cos nz = (-1)^m2^{n-1}\prod_{p=0}^{n-1}\ \cos(z+p\pi/n), 
$$
donc
$$
\cot nz = (-1)^m2^{n-1}\prod_{p=0}^{n-1}\ \cot(z+p\pi/n), 
$$
que l'on peut r\'e\'ecrire comme 
$$
\cot nz = (-1)^m2^{n-1}\prod_{p=-m}^{m}\ \cot(z-p\pi/n)\ .
\eqno{(3.3.1)} 
$$

{\bf 3.4.} On a~: 
$$
\cot(a+b) = \frac{\cos(a+b)}{\sin(a+b)} = \frac{
\cos a\cos b - \sin a\sin b}{\sin a\cos b + \cos a\sin b} 
$$
$$
= \frac{1 - \tan a\tan b}{\tan a + \tan b}\ .
$$
Donc
$$
\cot nz = (-1)^m2^{n-1}\prod_{p=-m}^{m}\ \frac{1 + \tan z\tan(p\pi/n)}
{\tan z - \tan(p\pi/n)}\ .
$$
Ceci est une fraction rationelle dont le num\'erateur est de degr\'e 
$n-1$ en $u = \tan z$ et le d\'enominateur est de degr\'e $n$, ayant 
les racines simples. Il s'en suit qu'on peut \'ecrire une 
d\'ecomposition en \'el\'ements simples~:
$$
\cot nz = \sum_{p=-m}^{m}\ \frac{a_p}{u - \tan(p\pi/n)}
$$
avec
$$
a_p = \lim_{z\ra p\pi/n}\cot nz\cdot(\tan z - \tan(p\pi/n)) 
$$
$$
= \lim_{z\ra p\pi/n} \frac{\cos nz}{\sin nz}\cdot 
\frac{\sin(z - p\pi/n)}{\cos z\cos(p\pi/n)} 
$$
$$
= \frac{1}{\cos^2(p\pi/n)}\lim_{h\ra 0}\frac{\cos(nh + p\pi)\sin h}
{\sin(nh + p\pi)} = 
\frac{1}{\cos^2(p\pi/n)}\lim_{h\ra 0}\frac{(-1)^p\sin h}
{(-1)^p\sin nh} 
$$
$$
= \frac{1}{n\cos^2(p\pi/n)}\ .
$$
Donc
$$
\cot nz = \sum_{p=-m}^{m}\ \frac{1}{n\cos^2(p\pi/n)(\tan z - \tan(p\pi/n))}\ . 
$$
En rempla\c{c}ant $z$ par $z/n$, 
$$
\cot z = \sum_{p=-m}^{m}\ \frac{1}{n\cos^2(p\pi/n)(\tan(z/n) - \tan(p\pi/n))} 
$$
$$
= \frac{1}{n\tan(z/n)} + \sum_{p=1}^m 
\frac{1}{n\cos^2(p\pi/n)}\cdot\frac{2\tan(z/n)}
{\tan^2(z/n) - \tan^2(p\pi/n)} 
$$
$$
= \frac{1}{n\tan(z/n)} + \sum_{p=1}^m\ \cdot\frac{2n\tan(z/n)}
{\cos^2(p\pi/n)(n\tan(z/n))^2 - (n\sin(p\pi/n))^2}\ .
$$
On a donc d\'emontr\'e le 

{\bf 3.5.} {\it Th\'eor\`eme.} Pour tout $n=2m+1$ impair 
$$
\cot z = \frac{1}{n\tan(z/n)} + \sum_{p=1}^m\ \cdot\frac{2n\tan(z/n)}
{\cos^2(p\pi/n)(n\tan(z/n))^2 - (n\sin(p\pi/n))^2}\ .
$$

En faisant $m\ra \infty$, on arrive \`a~: 

{\bf 3.6.} {\it Th\'eor\`eme.}
$$
\cot z = \frac{1}{z} + \sum_{p=1}^\infty\ \frac{2z}{z^2 - p^2\pi^2}\ .
$$

{\bf 3.7.} Revenons au d\'eveloppement de sinus 3.1. Supposons toujours que 
$n = 2m+1$ est impair. Alors 3.1 peut s'\'ecrire 
$$
\sin nz = (-1)^m 2^{n-1}\prod_{p=-m}^m\ 
\sin(z - p\pi/n) 
$$
$$
= (-1)^m 2^{n-1}\sin z\prod_{p=1}^m\ 
\sin(z - p\pi/n)\sin(z + p\pi/n)\ .
$$
On v\'erifie ais\'ement la formule suivante~:
$$
\sin^2(a+b) - \sin^2(a-b) = \sin 2a\sin 2b, 
$$
d'o\`u
$$
\sin a\sin b = \sin^2((a+b)/2) - \sin^2((a-b)/2)\ .
$$
Il s'en suit, 
$$
\sin(z - p\pi/n)\sin(z + p\pi/n) = \sin^2 z - \sin^2(p\pi/n), 
$$
d'o\`u
$$
\sin nz = 2^{n-1}\sin z\prod_{p=1}^m\ 
(\sin^2(p\pi/n) - \sin^2 z)\ .
$$
Or, d'apr\`es 3.2, 
$$
\prod_{p=1}^m\ \sin^2(p\pi/n) = \frac{n}{2^{n-1}}, 
$$
d'o\`u      
$$
\sin nz = n\sin z\prod_{p=1}^m\ 
(1 - (\sin^2z/\sin^2(p\pi/n)))\ .
$$
En rempla\c{c}ant $z$ par $z/n$, on arrive au 

{\bf 3.8} {\it Th\'eor\`eme.} Si $n=2m+1$ est impair alors
$$
\sin z = n\sin(z/n)\prod_{k=1}^m\ 
\biggl(1 - \frac{\sin^2(z/n)}{\sin^2(k\pi/n)}\biggr)\ .
$$

Maintenant si l'on fait tendre $m$ vers l'infini, on obtient le 

{\bf 3.9.} {\it Th\'eor\`eme.} 
$$
\sin z = z \cdot \prod_{p=1}^\infty \biggl(1 - \frac{z^2}{p^2\pi^2}\biggr)\ .
$$
(Convergence uniforme dans des sous-ensembles compacts.) 

\bigskip 

{\it Application aux nombres de Bernoulli} 

\bigskip 

{\bf 3.10.} On a~: 
$$
\cot(iz/2) = i\cdot\frac{e^{-z} + e^z}{e^{-z} - e^z} = 
- i\cdot \frac{e^z + 1}{e^z - 1}, 
$$
d'o\`u~: 
$$
\frac{z}{e^z - 1} = 
\frac{z}{2}\cdot \biggl(-1 + \frac{e^z + 1}{e^z - 1}\biggr) = 
- \frac{z}{2} + \frac{iz}{2}\cot(iz/2)\ .
$$
On rappelle que les nombres de Bernoulli sont d\'efinis par~: 
$$
\frac{z}{e^z - 1} = 1 - \frac{z}{2} + \sum_{n=1}^\infty\ 
b_{2n}\frac{z^{2n}}{(2n)!}\ .
$$

{\bf 3.11.} Le d\'eveloppement de $\cot$ nous dit~: 
$$
\cot z - \frac{1}{z} = \sum_{n=1}^\infty\ \frac{2z}{z^2 - n^2\pi^2}\ .
$$
Maintenant~: 
$$
\frac{2z}{z^2 - n^2\pi^2} = -\frac{2z}{n^2\pi^2}\cdot 
\frac{1}{1 - z^2/n^2\pi^2} = -\frac{2z}{n^2\pi^2}\cdot 
\sum_{k=0}^\infty\ \frac{z^{2k}}{n^{2k}\pi^{2k}} 
$$
$$
= - 2\sum_{k=1}^\infty\ \frac{z^{2k-1}}{n^{2k}\pi^{2k}}
$$
($|z| < \pi$). En \'echangeant l'ordre de sommations, 
il s'en suit~:
$$
\cot z = \frac{1}{z} - 2\sum_{k=1}^\infty\ \frac{S_{2k}}{\pi^{2k}} 
z^{2k-1},
$$ 
o\`u 
$$
S_k = \sum_{n=1}^\infty\ \frac{1}{n^k}\ .
$$
Donc 
$$
\frac{z}{e^z - 1} = - \frac{z}{2} + \frac{iz}{2}\cdot\biggl( 
\frac{2}{iz} + 2\sum_{k=1}^\infty\ \frac{S_{2k}}{\pi^{2k}} 
(-1)^ki\frac{z^{2k-1}}{2^{2k-1}}\biggr) 
$$
$$
= 1 - \frac{z}{2} + \sum_{k=1}^\infty\ 
(-1)^{k-1}\frac{S_{2k}}{2^{2k-1}\pi^{2k}}z^{2k}\ .
$$

{\bf 3.12.} En comparant avec 3.10, 
$$
b_{2n} = (-1)^{n-1}(2n)!\frac{2S_{2n}}{(2\pi)^{2n}}, 
$$
ou
$$
S_{2n} = (-1)^{n-1}\frac{(2\pi)^{2n}}{2(2n)!}b_{2n},  
$$
$n\geq 1$.  

\bigskip\bigskip

\centerline{\bf \S 4. Une formule de Ramanujan} 

\bigskip\bigskip

{\bf 4.1.} On agit \`a la Eisenstein.  
On suit [A], Chapitre II, no. 10. Commen\c{c}ons par le
d\'eveloppement de $\cot$~: 
$$
\pi \cot \pi u = \frac{1}{u} + \sum_{m\in\BZ, m\neq 0}\ \biggl(\frac{1}{u+m} -
\frac{1}{m}\biggr) 
$$
$$
= \frac{1}{u} + \sum_{m=1}^\infty\ \biggl(\frac{1}{u+m} + 
\frac{1}{u-m}\biggr)\ .
$$
On pose $w = e^{2\pi iu}$~; alors 
$$
\cot\pi u = i\frac{w+1}{w-1} = i\cdot\biggl(1 + \frac{2}{w-1}\biggr) = 
- i + 2i\sum_{n=1}^\infty\ w^n, 
$$
si $|w| < 1$, i.e. $\Im u > 0$. Il s'en suit, 
$$
\frac{1}{u} + \sum_{m=1}^\infty\ \biggl(\frac{1}{u+m} + 
\frac{1}{u-m}\biggr) = 
- \pi i - 2\pi i\sum_{n=1}^\infty\ w^n\ .
$$ 

{\bf 4.2.} On d\'erive $p$ fois par rapport \`a $u$~; 
puisque $(d/du)^p(w) = (2\pi i)^pw$, on a~: 
$$
(-1)^p p!\sum_{m=-\infty}^\infty\ \frac{1}{(u+m)^{p+1}} 
= - (2\pi i)^{p+1}\sum_{k=1}^\infty\ k^pw^k\ .
$$
On pose $u = n\tau$, $n > 0,\ \Im\tau > 0$,  
$$
(-1)^p p!\sum_{m=-\infty}^\infty\ \frac{1}{(m+n\tau)^{p+1}} 
= - (2\pi i)^{p+1}\sum_{k=1}^\infty\ k^pe^{2kn\pi i\tau}
$$
et l'on r\'ealise la somme sur $n$~:
$$
(-1)^p p!\sum_{n=1}^\infty\sum_{m=-\infty}^\infty\ \frac{1}{(m+n\tau)^{p+1}} 
= - (2\pi i)^{p+1}\sum_{k=1}^\infty\ k^p
\frac{e^{2k\pi i\tau}}{1 - e^{2k\pi i\tau}}
\eqno{(4.2.1)}
$$
(attention~: on a chang\'e l'ordre des sommations \`a droite.) 

{\bf 4.3.} Maintenant supposons que $p = 2l-1$ est impair et $p\geq 2$ 
(i.e. $l \geq 2$). On peut alors r\'e\'ecrire (4.2.1)~: 
$$
\frac{1}{2}\sum_{m,n}\ '\frac{1}{(m+n\tau)^{2l}} - 
\sum_{m=1}^\infty\ \frac{1}{m^{2l}} = \frac{(2\pi i)^{2l}}{(2l-1)!} 
\sum_{k=1}^\infty\ \frac{k^{2l-1}e^{2k\pi i\tau}}
{1 - e^{2k\pi i\tau}}\ .
$$ 
On utilise la notation  
$$
E_k(\tau) = \sum_{m,n}\ '\frac{1}{(m+n\tau)^{k}}
$$
pour les s\'eries d'Eisenstein.

{\bf 4.4.} Consid\'erons le cas sp\'ecial $\tau = i$~: 
$$
\frac{1}{2} E_{2l}(i) - \zeta(2l) = \frac{(-1)^l(2\pi)^{2l}}{(2l-1)!} 
\sum_{k=1}^\infty\ \frac{k^{2l-1}e^{-2k\pi}}{1 - e^{-2k\pi}} 
$$  
$$
= \frac{(-1)^l(2\pi)^{2l}}{(2l-1)!} 
\sum_{k=1}^\infty\ \frac{k^{2l-1}}{e^{2k\pi} - 1}\ .
$$
On rappelle en revanche que
$$
\zeta(2l) = (-1)^{l-1}\frac{(2\pi)^{2l}}{2(2l)!}b_{2l},  
$$
cf. 3.12. Il s'en suit~: 
$$
\sum_{k=1}^\infty\ \frac{k^{2l-1}}{e^{2k\pi} - 1} = 
(-1)^l\frac{(2l-1)!}{2(2\pi)^{2l}} E_{2l}(i) 
+ \frac{b_{2l}}{4l}\ .
$$

{\bf 4.5.} Supposons que $l = 2j+1$ est impair. Alors 
$$
E_{2l}(i) = \sum_{m,n}\ '\frac{1}{(m+ni)^{2l}} = (-i)^{2l} 
\sum_{m,n}\ '\frac{1}{(-mi+n)^{2l}} = - E_{2l}(i), 
$$    
donc $E_{2l}(i) = 0$. Il d\'ecoule que 
$$
\sum_{k=1}^\infty\ \frac{k^{2l-1}}{e^{2k\pi} - 1} = 
\frac{b_{2l}}{4l}
$$
dans ce cas. Ceci est une formule de Ramanujan. 

{\bf 4.6.} En g\'en\'eral, on d\'efinit les fonctions de Weierstrass~: 
$$
\sigma(u) = \sigma(\omega_1,\omega_2;u) = u\prod\ '\ 
\biggl(1 - \frac{u}{\omega}\biggr)e^{u/\omega + u^2/(2\omega^2)} \ ,
$$
o\`u $\omega = m\omega_1 + n\omega_2$ et
$$
\prod\ ' = \prod_{(m,n) \in \BZ^2 - \{(0,0)\} }  \ .
$$
Cette fonction est analogue de $\sin u$. Ensuite,  
$$
\zeta(u) = \frac{\sigma'(u)}{\sigma(u)} = 
\sum\ '\ \biggl\{\frac{1}{u - \omega} + \frac{1}{\omega} + 
\frac{u}{\omega^2}\biggr\} \ ,
$$
analogue de $\cot u$~; et
$$
\CP(u) = - \zeta'(u) = \sum\ '\ \biggl\{\frac{1}{(u-\omega)^2}
 - \frac{1}{\omega^2}\biggr\}\ ,
$$
analogue de $- \text{cosec}^2 u$. On a alors le d\'eveloppement 
de Laurent en $0$~:
$$
\zeta(u) = \frac{1}{u} - E_4u^3 - E_6u^5 - E_8u^7 - \ldots \ ,
$$
o\`u 
$$
E_n = E_n(\omega_1,\omega_2) = \sum\ '\ \frac{1}{\omega^n} \ .
$$
Donc 
$$
\CP(u) = \frac{1}{u^2} + 3E_4u^2 + 5E_6u^4 + 7E_8u^6 + \ldots 
$$
La fonction $\CP(u)$ satisfait les \'equations diff\'erentielles
$$
\CP^{\prime 2}(u) = 4\CP^3(u) - g_2\CP(u) - g_3\ ,
$$
$$
\CP''(u) = 6\CP^3(u) - g_3/2\ ,
$$
o\`u 
$$
g_2 = 60E_4,\ g_3 = 140E_6\ .
$$

{\bf 4.7.} Le cas du r\'eseau Gaussien $(\omega_1,\omega_2) = (1,i)$ 
a \'et\'e trait\'e par Hurwitz, [H]. On consid\`ere la fonction de Weierstrass 
qui satisfait l'\'equation diff\'erentielle 
$$
\CP^{\prime 2}(u) = 4\CP^3(u) - 4\CP(u), 
$$
donc $g_2 = 1$, $g_3 = 0$. On introduit la p\'eriode correspondante~: 
$$
\omega = 2\int_0^1\ \frac{dx}{\sqrt{1-x^4}}, 
$$
cf. une d\'efinition de $\pi$~: 
$$
\pi = 2\int_0^1\ \frac{dx}{\sqrt{1-x^2}} \ .
$$
On d\'efinit alors les nombres rationels $E_n$ par 
$$
\CP(u) = \frac{1}{u^2} + \frac{2^4E_1}{4}\cdot\frac{u^2}{2!} + 
\frac{2^8E_2}{8}\cdot\frac{u^6}{6!} + \ldots + 
\frac{2^{4n}E_n}{4n}\cdot\frac{u^{4n-2}}{(4n-2)!} + \ldots 
$$
On a $E_1 = 1/10$ et $E_n$ satisfait une relation de recurrence
$$
E_n = \frac{3}{(2n-3)(16n^2 - 1)}\ \sum_{k=1}^{n-1}\ 
(4k-1)(4n-4k-1)\binom{4n}{4k} E_kE_{n-k}\ .
$$
Alors 
$$
\sum\ '\ \frac{1}{(r+is)^{4n}} = \frac{(2\omega)^{4n}}{(4n)!}E_n\ .
$$  

\bigskip\bigskip

\centerline{\bf \S 5. Une int\'egrale de Legendre} 

\bigskip\bigskip

{\bf 5.1.} On rappelle que 
$$
\frac{t}{e^t - 1} = \frac{t}{2i}\cot\frac{t}{2i} - \frac{t}{2} = 
1 - \frac{t}{2} + \sum_{n=1}^\infty\ \frac{b_{2n}t^{2n}}{(2n)!} \ .
$$
On peut donc poser $b_{1} = - 1/2$. 

{\bf 5.2.} {\it Th\'eor\`eme} (Legendre). 
$$
\int_0^\infty\ \frac{\sin ax}{e^{2\pi x} - 1} dx = 
\frac{1}{4}\frac{e^a + 1}{e^a - 1} - \frac{1}{2a} \ .
$$

On donne deux d\'emonstrations. 

{\bf 5.3.} La {\it premi\`ere d\'emonstration} utilise le developpement de 
cotangent. 

On a~: 
$$
\frac{1}{e^{2\pi x} - 1} = e^{-2\pi x}\sum_{n=0}^\infty\ e^{-2\pi nx} = 
\sum_{n=1}^\infty\ e^{-2\pi nx}
$$
($x > 0$), d'o\`u
$$
I: = \int_0^\infty\ \frac{\sin ax}{e^{2\pi x} - 1} dx = 
\sum_{n=1}^\infty\ \int_0^\infty\ \sin ax\ e^{-2\pi nx} dx \ .
$$
Or,
$$
\int_0^\infty\ \sin ax\ e^{-2\pi nx} dx = 
\frac{1}{2i}\ \int_0^\infty\ (e^{iax} - e^{-iax})e^{-2\pi nx} dx\ ,
$$
o\`u
$$
\int_0^\infty\ e^{iax - 2\pi n x} dx = 
\frac{1}{ia - 2\pi n}e^{iax - 2\pi n x}\biggr|_0^\infty = 
\frac{1}{2\pi n - ia} = \frac{2\pi n + ia}{a^2 + 4\pi^2 n^2} \ .
$$
Donc
$$
\int_0^\infty\ \sin ax\ e^{-2\pi nx} dx = 
\frac{a}{a^2 + 4\pi^2 n^2}\ ,
$$
d'o\`u 
$$
I = \sum_{n=1}^\infty\  \frac{a}{a^2 + 4\pi^2 n^2}\ .
$$
Rappelons~: 
$$
\sum_{n=1}^\infty\ \frac{2a}{a^2 - \pi^2 n^2} = \cot a - \frac{1}{a}\ .
$$
De l\`a~: 
$$
\sum_{n=1}^\infty\  \frac{a}{a^2 + 4\pi^2 n^2} = 
\sum_{n=1}^\infty\  \frac{a/4}{a^2/4 + \pi^2 n^2} = 
\frac{1}{4i}\sum_{n=1}^\infty\  \frac{ia}{- (ia/2)^2 + \pi^2 n^2} 
$$
$$
= - \frac{1}{4i}\biggl(\cot(ia/2) - \frac{2}{ia}\biggr) = 
- \frac{1}{4i}\cot(ia/2) - \frac{1}{2a}\ .
$$
Or, 
$$
\cot(ia/2) = \frac{\cos(ia/2)}{\sin(ia/2)} = 
\frac{i(e^{-a/2} + e^{a/2})}{e^{-a/2} - e^{a/2}} = 
\frac{i(1 + e^{a})}{1 - e^{a}}\ ,
$$
donc
$$
- \frac{1}{4i}\cot(ia/2) = \frac{1}{4}\frac{e^a + 1}{e^a - 1}, 
$$
{\it quod erat demonstrandum}.

\bigskip

{\bf 5.4.} La {\it deuxi\`eme d\'emonstration} utilise la formule de Cauchy~; elle 
a \'et\'e propos\'ee comme un exercice dans [WW], Ch. 6, 6.4, Example 2. 
Le calcul a \'et\'e fait par Nabil Rachdi. 

On d\'efinit~:
$$
I_a(\epsilon,R) = \int_\epsilon^R\ \frac{e^{iax}}{e^{2\pi x} - 1}dx\ ;\ \ 
I_a(\epsilon) := I_a(\epsilon,\infty)\ .
$$
Alors
$$
I = \lim_{\epsilon\ra 0}\frac{I_a(\epsilon) - I_{-a}(\epsilon)}{2i}\ .
$$
Consid\'erons le contour "rectangulaire" $\Gamma = \Gamma(\epsilon,R)$ suivant~: 
$$
\Gamma = \bigcup_{i=1}^6\ \Gamma_i = \{\epsilon\leq z \leq R\}\cup 
\{z = R + it|\ 0\leq t \leq 1\} \cup 
\{z = t + i|\ R\geq t \geq \epsilon\} \cup
$$
$$
\cup \{ z = i + \epsilon e^{i\theta}| 0 \geq \theta \geq - \pi/2\}\cup 
\{z = it|\ 1 - \epsilon \geq t \geq \epsilon\} \cup 
\{ z = \epsilon e^{i\theta}|\ \pi/2 \geq \theta \geq 0\}\ .
$$
On pose~: 
$$
f(z) = \frac{e^{iaz}}{e^{2\pi z} - 1}\ .
$$
Puisque $e^{2\pi z} = 1$ ssi $z = ni,\ n\in\BZ$, cette fonction n'a pas de 
singularit\'es \`a l'int\'erieur de $\Gamma$, donc 
$$
0 = \int_\Gamma\ f(z) dz = \sum_{i=1}^6\ \int_{\Gamma_i}\ f(z)dz\ .
$$
Calculons les int\'egrales $\int_{\Gamma_i}\ f(z)dz$ s\'epar\'ement. On a~:
$$
\int_{\Gamma_1}\ f(z)dz = I(\epsilon,R)\ .
$$
De m\^eme,
$$
\int_{\Gamma_3}\ f(z)dz = - \int_\epsilon^R\ \frac{e^{ia(x+i)}}
{e^{2\pi(x+i)} - 1} dx = - e^{-a}I_a(\epsilon,R)\ .
$$
Ensuite, 
$$
\int_{\Gamma_2(R)}\ f(z)dz = \int_0^1\ \frac{e^{ia(R+it)}}{e^{2\pi(R+it)}-1} 
dt \lra 0\ \text{quand\ }R\ra\infty\ .
$$
Les int\'egrales sur les quarts de cercles~: 
$$
\int_{\Gamma_4(\epsilon)}\ f(z) dz = \int_0^{-\pi/2}\ 
\frac{e^{ia(i+\epsilon e^{i\theta})}}{e^{2\pi(i+\epsilon e^{i\theta})} 
- 1} i\epsilon e^{i\theta} d\theta 
$$
$$
= ie^{-a}\int_0^{-\pi/2}\ 
\frac{e^{ia\epsilon e^{i\theta}}\epsilon e^{i\theta}}
{e^{2\pi\epsilon e^{i\theta}} - 1} d\theta\ .
$$
Or, la fonction sous l'int\'egrale 
$$
\frac{e^{ia\epsilon e^{i\theta}}\epsilon e^{i\theta}}
{e^{2\pi\epsilon e^{i\theta}} - 1} \sim 
\frac{\epsilon e^{i\theta}}{2\pi\epsilon e^{i\theta}} = \frac{1}{2\pi}\ 
\text{quand\ }\epsilon \ra 0, 
$$
d'o\`u
$$
\int_{\Gamma_4(\epsilon)}\ f(z) dz \lra - \frac{i}{4}e^{-a}\      
\text{quand\ }\epsilon \ra 0\ .
$$
De m\^eme, pour $\int_{\Gamma_6}$ on trouve~: 
$$
\int_{\Gamma_6(\epsilon)}\ f(z) dz \lra - \frac{i}{4}\      
\text{quand\ }\epsilon \ra 0\ .
$$

{\bf 5.5.} Finalement, il reste \`a traiter l'int\'egrale 
$\int_{\Gamma_5}$. On a~: 
$$
\int_{\Gamma_5(\epsilon)}\ f(z)dz = - \int_\epsilon^{1-\epsilon}\ 
\frac{e^{ia\cdot it}}{e^{2\pi it} - 1} i dt = 
- i\int_\epsilon^{1-\epsilon}\ 
\frac{e^{-at}}{e^{2\pi it} - 1} dt := J_a(\epsilon)\ .
$$
Par la formule de Cauchy, 
$$
0 = (1 - e^{-a})I_a(\epsilon) - \frac{i}{4}(1 + e^{-a}) + J_a(\epsilon) 
+ o(\epsilon),
$$
d'o\`u, en posant $y = e^a$, 
$$
I_a(\epsilon) = \frac{i}{4}\cdot \frac{1 + y^{-1}}{1 - y^{-1}} - 
\frac{1}{1 - y^{-1}} J_a(\epsilon) + o(\epsilon)\ .
$$
Il s'en suit~:
$$
I_{-a}(\epsilon) = \frac{i}{4}\cdot \frac{1 + y}{1 - y} - 
\frac{1}{1 - y} J_{-a}(\epsilon) + o(\epsilon)\ .
$$
Or,
$$
J_{-a}(\epsilon) = - i\int_\epsilon^{1-\epsilon} 
\frac{e^{at}}{e^{2\pi it} - 1} dt 
$$
($x = - t + 1$)
$$
= - i\int_{1-\epsilon}^{\epsilon} 
\frac{e^ae^{-ax}}{e^{-2\pi ix} - 1}\cdot (- dx) = - iy
\int_\epsilon^{1-\epsilon} \frac{e^{2\pi ix}\cdot e^{-ax}}
{1 - e^{2\pi ix}} dx 
$$
$$
= -iy \int_\epsilon^{1-\epsilon} e^{-ax}\cdot \biggl(-1 + 
\frac{1}{1 - e^{2\pi ix}}\biggr) dx = yJ_a(\epsilon) + 
iy\int_\epsilon^{1-\epsilon} e^{-ax} dx \ .
$$
La derni\`ere int\'egrale
$$
\int_\epsilon^{1-\epsilon} e^{-ax} dx = \int_0^1 e^{-ax} dx + o(\epsilon) = 
- \frac{e^{-ax}}{a} \biggr|_0^{1} + o(\epsilon) = \frac{1-y^{-1}}{a} 
+ o(\epsilon)\ .
$$
Ainsi, 
$$
J_{-a}(\epsilon) = - yJ_a(\epsilon) + \frac{i(y-1)}{a} + o(\epsilon)\ .
\eqno{(5.5.1)}
$$
On obtient~: 
$$
I_a(\epsilon) -  I_{-a}(\epsilon) = \frac{i}{4}\cdot \biggl[
\frac{1 + y^{-1}}{1 - y^{-1}} - \frac{1 + y}{1 - y}\biggr] 
+ \frac{J_{-a}(\epsilon)}{1-y} - \frac{J_a(\epsilon)}{1-y^{-1}} 
+ o(\epsilon) \ .
$$
Ici~:
$$
\frac{1 + y^{-1}}{1 - y^{-1}} - \frac{1 + y}{1 - y} = 
2 \cdot \frac{y+1}{y-1} 
$$
et
$$
\frac{J_{-a}(\epsilon)}{1-y} - \frac{J_a(\epsilon)}{1-y^{-1}} = 
- \frac{yJ_a(\epsilon)}{1-y} - \frac{i}{a} 
- \frac{yJ_a(\epsilon)}{y-1} + o(\epsilon) = 
- \frac{i}{a} + o(\epsilon)\ .
\eqno{(5.5.2)}  
$$
On peut voir en (5.5.1), ou en la formule \'equivalente (5.5.2), 
une \'equation fonctionelle pour la fonction $J_a(\epsilon)$~; 
remarquons que l'int\'egrale $J_a(\epsilon)$ diverge quand 
$\epsilon\ra 0$. En revenant \`a $I_a$, on obtient~: 
$$
I_a(\epsilon) -  I_{-a}(\epsilon) = \frac{i}{2}\cdot \frac{y+1}{y-1} 
- \frac{i}{a} + o(\epsilon), 
$$
d'o\`u
$$
\frac{I_a(\epsilon) -  I_{-a}(\epsilon)}{2i} = 
\frac{1}{4}\cdot \frac{y+1}{y-1} 
- \frac{1}{2a} + o(\epsilon)\ .
$$
En faisant tendre $\epsilon$ vers z\'ero, on obtient la valeur 
de l'int\'egrale de Legendre. 

\bigskip 
 
{\bf 5.6.} {\it Th\'eor\`eme.} Pour $n\geq 1$, 
$$
(-1)^{n-1}b_{2n} = 4n \int_0^\infty\ \frac{t^{2n-1}}{e^{2\pi t} - 1} dt\ .
\eqno{(5.6.1)}
$$

On peut consid\'erer cela comme une {\it deuxi\`eme d\'efinition} 
des nombres de Bernoulli (Jacob Bernoulli, {\it Ars conjectandi}, 1713, p. 97). 

On en donne deux d\'emonstrations. 

La {\it premi\`ere d\'emonstration} utilise les valeurs de $\zeta(s)$ en 
points positifs pairs (donc le d\'eveloppement de $\cot$)~: 
$$ 
\int_0^\infty\ \frac{t^{2n-1}}{e^{2\pi t} - 1} dt = 
\int_0^\infty\ t^{2n-1}e^{-2\pi t}\sum_{k=0}^\infty e^{-2\pi kt} dt 
$$
$$
= \int_0^\infty\ t^{2n-1}\sum_{k=1}^\infty e^{-2\pi kt} dt = 
\sum_{k=1}^\infty \int_0^\infty\ t^{2n-1} e^{-2\pi kt} dt 
$$
($x = 2\pi kt$)
$$
= \sum_{k=1}^\infty \int_0^\infty\ (x/(2\pi k))^{2n-1} e^{-x} dx/(2\pi k) 
$$
$$
= (2\pi)^{-n}\Gamma(2n) \sum_{k=1}^\infty\ \frac{1}{k^{2n}} = 
(2\pi)^{-n}(2n-1)! S_{2n} 
$$
(voir 3.12)
$$
= (-1)^{n-1}\frac{b_{2n}}{4n}\ .
$$
   
\bigskip

{\bf 5.7.} La {\it deuxi\`eme d\'emonstration} utilise l'int\'egrale de Legendre (avec 
la preuve par la formule de Cauchy), cf. [WW], 7.2~: 
$$  
\int_0^\infty\ \frac{\sin ax}{e^{\pi x} - 1} dx = 
- \frac{1}{2a} + \frac{i}{2}\cot ia = \frac{1}{2a} 
\sum_{n=1}^\infty\ b_{2n} \frac{(2a)^{2n}}{(2n)!}\ .
$$
En d\'erivant $2n$ fois et en posant $a=0$ et $x=2t$, on en d\'eduit (5.6.1). 

\bigskip 
 
En particulier, si $n$ est impair, 
$$
\frac{b_{2n}}{4n} = \int_0^\infty\ \frac{t^{2n-1}}{e^{2\pi t} - 1} dt,
$$
cf. 4.5. On arriv\'e ainsi \`a l'assertion~:  

{\bf 5.8. Theorema pulcherissimum.} {\it Si $n > 1$ est un entier impair, alors 
$$
\int_0^\infty\ \frac{t^{2n-1}}{e^{2\pi t} - 1} dt = 
\sum_{k=1}^\infty \frac{k^{2n-1}}{e^{2\pi k} - 1} = 
\frac{b_{2n}}{4n}\ .
$$}

\bigskip\bigskip

\newpage

\centerline{\bf Bibliographie}

\bigskip\bigskip

[A] N.I.Akhiezer, \'Elements de la th\'eorie des fonctions elliptiques 
(en russe). 2-\`eme \'edition, Nauka, Moscou, 1976.  

[B] N.Bourbaki, Fonctions d'une variable r\'eelle. Hermann, 1961.  

[D] R.Dedekind, Erl\"auterungen zu den Fragmenten XXVIII.
Dans~: Bernhard Riemann, 
Gesammelte Mathematische Werke, Herausgegeben unter Mit\-wirkung von Richard 
Dedekind und Heinrich Weber, Teubner, Leipzig, 1892, pp. 466 - 478.

[H] A.Hurwitz, \"Uber die Entwicklungscoefficienten der lemniscatischen 
Functionen. {\it Nachrichten von der K\"onigl. Gesellschaft der Wissenschaften 
zu G\"ottingen, Mathematisch - physikalische klasse}, 1897, 273 - 276. 

[Ram] S.Ramanujan, Modular equations and approximations to $\pi$.
{\it Quarterly Journal of Mathematics}, {\bf 14} (1914), pp. 350 - 372. 
(Dans~: Collected Papers of Srinivasa Ramanujan, AMS Chelsea Publishing, 2000, 
pp. 23 - 39.)   

[R] B.Riemann, Fragmente \"uber die Grenzfalle der elliptischen 
Modulfunctionen. Werke (Leipzig, 1892), pp. 455 - 465. 

[Sch] O.Schl\"omilch, \"Uber einige unendliche Reihen. {\it Ber. Verh. K. 
Sachs. Gesell. Wiss.} Leipzig, {\bf 29} (1877), pp. 101 - 105. 

[S] C.L.Siegel, A simple proof of $\eta(-1/\tau) = \eta(\tau)\sqrt{\tau/i}$.
{\it Mathematika}, {\bf 1} (1954), p. 4. 

[WW] E.T.Whittaker, G.N.Watson, A course of modern analysis. Fourth Edition, 
Cambridge University Press, 1927.

\bigskip

\enddocument